\documentclass[a4paper,11pt]{article}

\usepackage{amsmath,amssymb,amsfonts}
\usepackage{graphicx}
\usepackage{textcomp}
\usepackage{xcolor}

\usepackage[utf8]{inputenc}
\usepackage[T1]{fontenc}
\usepackage{lmodern}

\usepackage{hyperref}
\hypersetup{
    colorlinks,
    linkcolor={blue!55!black},
    citecolor={blue!50!black},
    urlcolor={blue!80!black}
}

\usepackage{tikz}
\usetikzlibrary{shapes.symbols}
\usetikzlibrary{decorations.pathreplacing}
\usepackage{amsthm}

\usepackage[margin=2cm]{geometry}

\newtheorem{prop}{Proposition}

\usepackage{mdframed}
\usepackage[font=sf]{caption}

\newcommand{\dps}{\displaystyle}
\newcommand{\E}[1]{\mathbb{E}[ #1 ]}
\newcommand{\bE}[1]{\mathbb{E}\left[ #1 \right]}

\newcommand{\sn}{k}           
\newcommand{\bl}{u}           
\newcommand{\blv}{\mathbf{u}} 
\newcommand{\st}{z}           
\newcommand{\stv}{\mathbf{z}} 
\newcommand{\cf}{c}           

\newcommand{\vf}{v}           
\newcommand{\wf}{w}           

\newcommand{\gs}{m}           
\newcommand{\dt}{\delta}      
\newcommand{\ev}{\mathbf{e}}  

\newcommand{\CV}{C}           
\newcommand{\jsn}{k}

\begin{document}

\title{{\Large\bf\sffamily Technical Report}\\[5pt]\sffamily
    Optimal Size-Aware Dispatching Rules via Value Iteration\\
  and Some Numerical Investigations}
\author{\sffamily Esa Hyytiä\footnote{Department of CS, University of Iceland, esa@hi.is}{ }
  and Rhonda Righter\footnote{Department of IEOR, UC Berkeley, rrighter@ieor.berkeley.edu}}
\date{\sffamily\small\today}
\maketitle

\begin{abstract}
  This technical report explains how optimal size-aware dispatching
  policies can be determined numerically using value iteration. It also contains
  some numerical examples that shed light on the nature of the optimal policies themselves.
  The report complements our {\em ``Towards the Optimal Dynamic Size-aware Dispatching''} article
  that will appear in Elsevier's Performance Evaluation in 2024 \cite{hyytia-peva-2024}.
\end{abstract}

\begin{mdframed}[backgroundcolor=yellow!10]
  \small\vspace{-5mm}
  \tableofcontents
\end{mdframed}

\section{Introduction}

\begin{figure}[t]
  \centering
  \begin{tikzpicture}[font=\sffamily\scriptsize]
    \draw[line width=1pt,->,>=latex]
    (0,0) node [left] {new jobs} --  (1,0)
    node [right,draw,thick,circle,outer sep=2pt,inner sep=2pt,fill=yellow!10,label=below:dispatcher\;]
    (dp) {$\alpha$};
    \foreach \i/\u in {1/2,2/2.6,3/1.5} {
      \fill[fill=blue!5] (3,\i/2-1-0.2) -- ++(1,0) --++(0,0.4) --++ (-1,0);
      \draw[thin,fill=brown!35] (4,\i/2-1-0.2) rectangle ++(-\u/4,0.4);
      \draw[thick] (3,\i/2-1-0.2) -- ++(1,0) --++(0,0.4) --++ (-1,0);
      \node at (3,\i/2-1) (s\i) {};
      \draw[fill=yellow!10] (4.25,\i/2-1) circle (0.2);
    }
    \draw[shorten >=2pt,>=latex,->,line width=1pt] (dp) -- ++ (0.9,0.5)  -- (s3);
    \draw[shorten >=2pt,>=latex,->,line width=1pt] (dp)                  -- (s2);
    \draw[shorten >=2pt,>=latex,->,line width=1pt] (dp) -- ++ (0.9,-0.5) -- (s1);
    \node at (3.5,1.0) {servers};
  \end{tikzpicture}
  \caption{Dispatching system.}
  \label{fig:dispatching-sys}
\end{figure}
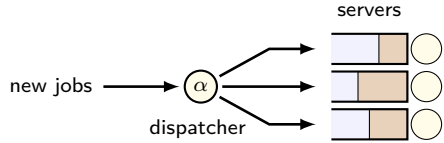

Parallel server systems where jobs are assigned to servers immediately upon arrival
are known as dispatching systems.
In the context of computing,
this is often referred to as load balancing (cf.\ multi-core CPUs, large data centers and cloud computing).
Dispatching is important also in manufacturing and service systems
(e.g., hospitals, help desks, call centers),
and various communication systems (cf.\ packet and flow level routing).
An example dispatching system is depicted in Figure~\ref{fig:dispatching-sys}.
It consists of three main elements:
i) jobs arrive according to some arrival process to dispatcher,
ii) dispatcher forwards each job immediately to a suitable server
utilizing information on job, state of the servers
and the objective function if available, and
iii) servers process jobs according to local scheduling rule.
In this paper, we consider a general arrival process,
size- and state-aware system (i.e., dispatching decisions can be based on
full state information about the job size and present backlogs in the available servers),
and servers process jobs in the first-come-first-served (FCFS) order.

This report consists of two parts.
First we discuss how the value iteration can be implemented,
which involves state-space truncation and discretization steps.
Then we discuss some numerical experiments that shed light on the nature
of the optimal dispatching rules, and in part
address some questions raised in \cite{hyytia-questa-2022}.

\subsection{Model for dispatching system}
The model for the dispatching system is depicted in Figure~\ref{fig:dispatching-sys}
and can be summarized as follows.
We have $\sn$ parallel FCFS servers and a dispatcher that assigns incoming
jobs to them immediately upon arrival. In the size-aware setting,
this decision may depend on
the current backlog $\bl_i$ in each server $i$
as well as the size of the current job.
For example the individually optimal dispatching policy least-work-left (LWL)
(see, e.g., \cite{foss-siberian-1980,harchol-balter-pdc-1999,akgun-aap-2013})
simply chooses the server with the shortest backlog thus
minimizing the latency for the present job.
Other popular heuristic dispatching policies
include JSQ \cite{haight-bio-1958,ephremides-tac-1980,whitt-or-1986}
and Round-Robin \cite{arian-orl-1992,down-queueing-2006,wu-ejor-2009,hyytia-itc32-2020-star,anselmi-tpds-2020}.

In this report, we assume that the job arrival process is defined
by a sequence of i.i.d.\ inter-arrival times $A_i \sim A$.
Similarly, the job sizes $X_i$ are i.i.d., $X_i \sim X$.
Costs are defined by an arbitrary cost function $\cf(\bl_i,x)$, where
$\bl_i$ denotes the backlog in server $i$ and $x$ is the size of the job.
The objective is to assign jobs to servers so as to minimize long-run average costs.

In our examples, we assume exponential distribution for both $A$ and $X$, i.e.,,
jobs arrive according to a Poisson process with rate $\lambda$ and
job sizes obey exponential distribution.
We let $f(x) = \mu e^{-\mu x}$ denote the pdf of job-size distribution,
and $f_A(t) = \lambda \,e^{-\lambda t}$ denotes the pdf of inter-arrival times.
For simplicity, we set $\mu=1$, i.e., $X \sim \text{Exp}(1)$.
The objective is to minimize the mean waiting time,
and thus our cost function is $\cf(\bl_i,x)=\bl_i$.

\subsection{Value functions and value iteration}
The optimal size-aware policy can be described by the value function $\vf(\bl_1,\ldots,\bl_{\sn})$.
In general, value functions are defined as expected deviations from the mean costs.
In our context, the stochastic process evolves in continuous time
and it is convenient to define value functions for two different time instants.

Let $\alpha$ denote an arbitrary dispatching policy.
The so-called post-assignment value function is defined for the time instant
when a new job has just been assigned to some server. Formally,
\begin{equation}\label{eq:vf-ggk-v}
\vf^{\alpha}(\blv) = \bE{ \sum_{i=1}^\infty\left(  \CV^{\alpha}_i(\blv) - \E{\CV^{\alpha}}\right) } + {c}_0,
\end{equation}
where the $\CV^{\alpha}_i(\blv)$'s
denote the costs the jobs arriving in the future incur with policy $\alpha$,
$\E{\CV^{\alpha}}$ is the corresponding mean cost incurred,
and ${c}_0$ is an arbitrary constant.

The pre-assignment value function is defined for the time instant when
a new job has just arrived, but it has not been assigned yet. Moreover,
its size is still unknown (while the immediately following dispatching decision
will be aware of it).
Formally,
\begin{equation}
  \label{eq:vf-ggk-w}
  \wf^{\alpha}(\blv) := \bE{ \sum_{i=0}^\infty \left(\tilde{\CV}^{\alpha}_i(\blv) - \E{\CV^{\alpha}} \right) } + \tilde{c}_0,
\end{equation}
where $\tilde{\CV}^{\alpha}_i(\blv)$
denotes the cost incurred by the $i$th job when the first job (arrival $0$)
will be assigned immediately in state $\blv$ (and change the state),
and $\tilde{c}_0$ is an arbitrary constant.

For optimal policy $\alpha^*$ we drop $\alpha$ from the superscripts, and write, e.g.,
\begin{equation}\label{eq:vf-opt}
  \vf(\blv) = \bE{ \sum_{i=1}^{\infty} \left( {\CV_i}(\blv) - \E{\CV} \right) } + c_0.
\end{equation}
The value functions corresponding to the optimal policy satisfy Bellman equation \cite{hyytia-peva-2024}:
\begin{prop}[Bellman equation]
  The pre-assignment value function $\wf(\blv)$ of the optimal dispatching policy satisfies
  \begin{equation}\label{eq:dynprog-ggk3}
    \wf(\blv) = \E{ \min_i \left\{ \cf(\bl_i,X) - \E{\CV} + \vf(\blv+\ev_i X) \right\} },
  \end{equation}
  where $\vf(\blv)$ is the post-assignment value function, for which it holds that,
  \begin{equation}\label{eq:vf-wz2}
    \vf(\blv) = \E{ \wf( (\blv-A \ev)^+) },
  \end{equation}
  where
  $(\cdot)^+$ denotes the component-wise maximum of zero and the argument.
\end{prop}

\begin{figure}
  \centering
  \begin{tikzpicture}[scale=0.3,font=\scriptsize,node distance=8mm]
    \draw[<->,>=latex,thick]
    (-9,5.5)
    -- node [rotate=90,above] {backlog in server 2, $\bl_2$}
    ++(0,-12.0) 
    -- node [below] {backlog in server 1, $\bl_1$}
    ++(18,0) ;
    %
    \begin{scope}[shift={(-9,-6)}]
    \node[draw,circle,minimum width=4pt,inner sep=0pt,fill,label=left:$\blv$] at (4,3) (a0) {};
    \node[minimum width=4pt,inner sep=0pt] at (4,2) (a01) {};
    \draw (4,8.5) node (tg1) {};
    \draw (9.5,3) node (tg2) {};
    \draw[thick, ->,>=latex,shorten >=3pt] (a0) -- (tg1);
    \draw[thick, ->,>=latex,shorten >=3pt] (a0) -- (tg2);
    \node[circle,minimum width=4pt,inner sep=0pt,draw,fill=red!50] at (4,6) (b1) {};
    \node[circle,minimum width=4pt,inner sep=0pt,draw,fill=red!50] at (7,3) (b2) {};
    \draw[line width=3.0pt,color=red!50!black,opacity=0.36] (b1) -- (a0) -- (b2);
    \node[xshift=15mm,yshift=6mm,anchor=west,inner sep=0pt] at (b1) (b1l) {$\cf(\bl_2,x)-\E{\CV} + \vf(\blv+x \ev_2)$};
    \draw[<-,shorten >=6pt,shorten <=6pt,thin] (b1) -- (b1l.south west);
    \node[below of=b1l,yshift=1mm,anchor=south,inner sep=1pt] (b2l) {
      $\cf(\bl_1,x)-\E{\CV} + \vf(\blv+x \ev_1)$};
    \draw[<-,shorten >=4pt,shorten <=6pt,thin] (b2) -- (b2l.south west);
    \draw [decorate,decoration={brace,raise=2pt}] (a0.north east) -- node[above,yshift=5pt] {$x$} (b2.north east);
    \draw [decorate,decoration={brace,raise=2pt}] (a0.north west) -- node[left,xshift=-5pt] {$x$} (b1.west);
    \end{scope}
  \end{tikzpicture}
  \quad
  \begin{tikzpicture}[scale=0.3,font=\scriptsize]
    \shade[
      left color = yellow!40,
      right color = red!40,
      shading angle = 135,
      opacity=90
    ]   (0,0) -- ++(4.5,0) --++(-5,-5) -- ++ (-4.5,0) --++(0,4.5) --++(5,5) --++(0,-4.5);
    \draw[<->,>=latex,thick]
    (-9,5.5)
    -- node [rotate=90,above] {backlog in server 2, $\bl_2$}
    ++(0,-12.0) 
    -- node [below] {backlog in server 1, $\bl_1$}
    ++(18,0) ;
    
    \draw[->,>=latex,ultra thick] (0,0) -- ++(5,0);
    \draw[->,>=latex,ultra thick] (0,0) -- ++(0,5);
    \draw[decorate,decoration={brace,raise=4pt}]  (0.3,0) -- node [above,xshift=8pt,yshift=5pt] {Job size $X$} ++(4.5,0);
    \foreach \i in {1,2,3,4} {
      \draw[->,>=latex] (\i,0) --++(-4.5,-4.5);
      \draw[->,>=latex] (0,\i) --++(-4.5,-4.5);
    }
    \draw[->,>=latex] (0,0) --++(-4.5,-4.5);

    \draw[decorate,decoration={brace,raise=5pt}] (4.5,0) -- node [below right,xshift=4pt,yshift=-4pt] {IAT $A$} ++(-4.8,-4.8);

    \node[draw,fill=black,circle,inner sep=1.5pt,minimum width=0mm]  at (0,0) (z0) {};
    \draw[<-,shorten <=7pt] (z0) --++(1.1,2.3) node [above right] {State $\blv$};
  \end{tikzpicture}

  \caption{One decision epoch consists of the job assignment of length
    $X$
    and the following IAT
    $A$.
    In the value iteration, the new value of state $\blv$ depends on the current values at the shaded area.}
\end{figure}
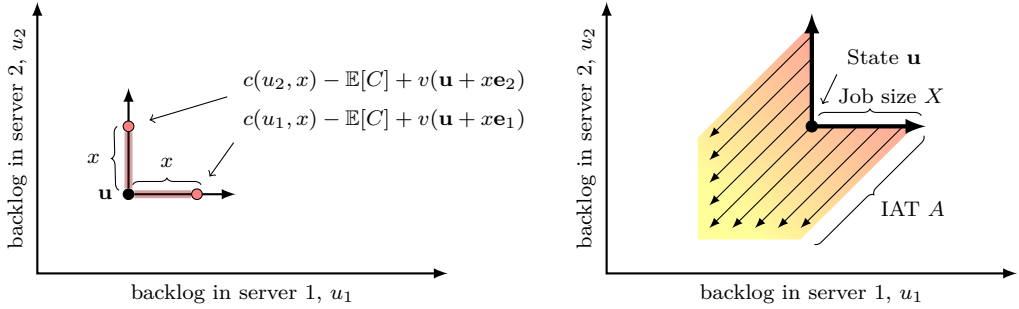

Note that \eqref{eq:dynprog-ggk3}
follows from requirement that the present job (of unknown size $X$) must be assigned in the optimal manner
minimizing the sum of the immediate cost, $\cf(\bl_i,X) - \E{\CV}$,
and the expected future costs, $\vf(\blv+\ev_i X)$.
Similarly, \eqref{eq:vf-wz2} follows from conditioning on the next job arriving after time $A$
and then utilizing the definition of $\wf(\blv)$.
The role of $(\cdot)^+$ operator is to ensure that no backlog becomes negative.

Consequently, the optimal value function can be obtained via the process known as value iteration.
Let $\vf_j(\blv)$ denote the value function at round $j$. Then, from \cite{hyytia-peva-2024}, we have:
\begin{equation}\label{eq:g-viter-b}
  \begin{array}{r@{\,}l}
    w_0     &= \dps \int\limits_0^\infty  f(x) \vf_{j}(\ev_1 x)  \, dx,\\
    \wf(\blv) &= \dps \int\limits_0^\infty  f(x) \min_i \left\{ \cf(\bl_i,x) + \vf_{j}(\blv+\ev_i x) \right\} \, dx - w_0,\\
    \vf_{j+1}(\blv) &= \dps \int\limits_0^\infty f_A(t) \wf(\blv-t \ev)\,dt,
  \end{array}
\end{equation}
where $w_0$ corresponds to the mean cost per job, $w_0 = \E{ \CV }$.

Each iteration round gives a new value function $\vf_{j+1}(\blv)$ based on $\vf_j(\blv)$.
Moreover, if the initial values correspond to some heuristic policy,
$w_0$ should decrease until the optimal policy has been found.
Numerical computations using a discrete grid and numerical integration,
however, introduce inaccuracies
that the stopping rule for value iteration should take into account.


\section{Implementing the value iteration}
\label{sec:time-step}

In this section,
we study how value iteration can be implemented using numerical methods.
Numerical \\
\parbox{105mm}{
computations should be stable and ideally the iteration
should converge fast.
Often when problems in continuous space are solved numerically one resorts
to discretization. Also our computations are carried out using a {\em finite grid} of size $\gs^\sn$,
where $\gs$ defines the length of each dimension and $\sn$ is the number of servers.
Hence, both $\vf[\stv]$ and $\wf[\stv]$ are $\sn$-dimensional arrays
having $\gs^\sn$ elements. This is illustrated in
the figure on the right for $\sn=3$ servers.
}
\hfill
\parbox{6cm}{
  {
  \centering
  \includegraphics[width=4cm]{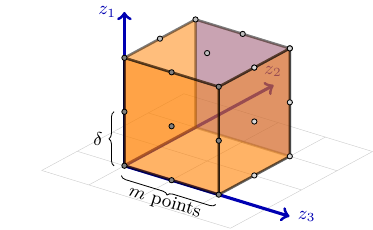}
  }\\
  \sffamily\small\parskip0pt
  Figure: Discretization and truncation of the state space yields $\sn$-dimensional array.
}

\paragraph{Discretization step $\boldsymbol{\dt}$:}
Let $\dt$ denote the discretization step or time step
for the grid. It is generally chosen to match the job size distribution $X$.
Basically, the discretization should be such that numerical integration, e.g., using
Simpson's method, gives reasonably accurate results.
That is, $\dt$ should be small enough so that
$$
\E{X} \approx \frac{\dt}{3} \left( h_0 + 4\,h_1 + 2\,h_2 + \ldots + h_{\jsn-1} \right),
\qquad
\text{where $h_i = i\dt \cdot f( i\,\dt )$},
$$
but not much larger, as otherwise the computational burden increases unnecessarily.
With Exp$(1)$-distribution,
a computationally efficient compromise is, e.g,  $\dt=0.25$.

The choice of a suitable discretization step depends also on the inter-arrival distribution $A$.
We will discuss this later in Section~\ref{sec:recursive-high-lambda}.

\paragraph{Value iteration algorithm:}
Let $\ev_i$ denote
a unit vector with the $i$th element equal to one and the rest zero.
For example, with $\sn=2$ servers, $\ev_1=(1,0)$ and $\ev_2=(0,1)$.
Similarly, without any subscript, $\ev$
is the all-one vector, $\ev=(1,\ldots,1)$.
The basic value iteration algorithm is given in Table~\ref{tbl:val-iter}.
Note that we have ignored the boundaries.
The crude way to deal with them is to restrict the index to the arrays
by defining for each coordinate $\stv_i$ a truncation
to $\{0,\ldots,\gs-1\}$:
$$
\tilde{\stv}_i = \max \{ 0,\, \min\{\gs-1, \stv_i\} \},
$$
and use $\tilde{\stv}$ instead of $\stv$.
Moreover,
instead of integrating from zero to
infinity, we simply truncate at some high value
where the integrand is very small.

\begin{table}
\begin{mdframed}[backgroundcolor=yellow!10]
\begin{itemize}
  \item {\bf Initialization:} Initially the value function can be set according
    to the random split (RND), for which the value function is a quadratic polynomial
    (see, e.g., \cite{hyytia-ejor-2012,hyytia-peva-2020}),
    $$ \vf[\stv] =  \sum_i \frac{ \lambda' (\stv_i \dt)^2}{2(1-\lambda')},$$
    where $\lambda'= \lambda/\sn$  is the server-specific arrival rate (uniform random split).
    Alternatively, we can initialize $\vf[\stv]$ simply with zeroes (or some suitable random values),
    $$
    \vf[\stv] =  0, \qquad \forall\, \stv.
    $$
    
  \item {\bf Iteration round:} One iteration round computes new values $\tilde{\vf}[\stv]$
    from a given $\vf[\stv]$:
    \begin{enumerate}
    \item
      First we apply Simpson's rule to determine $w_0$
      $$
      w_0 = \frac{\dt}{3} \left( h_1 + 4 h_2 + 2 h_3 + 4 h_4 + \ldots \right),
      $$
      where
      $$
      h_j = f(j\dt) \cdot \left( \cf(0,j\dt) + \vf[ j \ev_1] \right).
      $$
    \item
      Similarly,
      $\wf[\stv]$ can be updated by evaluating numerically the corresponding integral,
      $$
      \wf[\stv] =
      \frac{\dt}{3} \left( h_1 + 4 h_2 + 2 h_3 + 4 h_4 + \ldots \right),
      \qquad \forall\, \stv,
      $$
      where 
      $$
      h_j
      = f(j\dt) \cdot \left(
      \min_i \left\{ \cf(\st_i\dt,j\dt) + \vf[\stv+j\ev_i] \right\}
      - w_0 \right).
      $$
    \item
      Once $\wf[\stv]$ has been updated, we compute the new values
      for the $\vf[\stv]$, similarly by evaluating the
      corresponding integrals numerically, e.g., utilizing the same Simpson's composite rule,
      $$
      \tilde{v}(\stv) = \frac{\dt}{3} \left[
        h_1 + 4 h_2 + 2 h_3 + 4 h_4 + \ldots \right],
      $$
      where, assuming Poisson arrival process with rate $\lambda$,
      $$
      h_j
      = f_A(j\dt) \cdot \wf[ \stv - j \ev]
      = \lambda e^{-\lambda j\dt} \cdot \wf[ (\stv - j \ev)^+].
      $$
    \end{enumerate}
  \item{\bf Stopping rule:}
    The value iteration is repeated until the values in $\vf[\stv]$ converge.
    In practice, we can proceed, e.g., as follows:
    \begin{enumerate}
    \item Iterate first say 10-100 rounds to get started.
    \item Then iterate while $w_0$ changes (or decreases)
    \end{enumerate}
    The stopping criteria for iteration can be based also on
    the average change in the $\vf[\stv]$, which captures better the
    state of the value iteration than $w_0$, which reflects the change in
    one boundary of the whole state space.
    We return to this in our numerical examples.
\end{itemize}
\end{mdframed}
\caption{The basic value iteration algorithm.}
\label{tbl:val-iter}
\end{table}

\subsection{Recursive algorithm for Poisson arrivals}
\label{sec:recursive-poisson}

The Poisson arrival process is an important special case.
As mentioned in \cite{hyytia-peva-2024}, we can exploit the lack of memory
property of the Poisson process when updating $\vf[\stv]$.
Namely, 
$$
  \vf_{j+1}(\blv)
  =
  \int_0^h \lambda\,e^{-\lambda t}\,\wf(\blv-t\ev)\,dt
  +
  \int_h^\infty \lambda\,e^{-\lambda t}\,\wf(\blv-t\ev)\,dt,
$$
where the latter integral is equal to
$$
e^{-\lambda h}\cdot \vf_{j+1}(\blv-h\ev).
$$
Let $A(\blv)$ denote the first integral for interval $(0,h)$
\begin{equation}\label{eq:Au}
A(\blv) :=
\int_0^{h}
\lambda\,e^{-\lambda t}\,\wf(\blv-t\ev)\,dt.
\end{equation}
At grid point $\stv$, with $h=\dt$, 
$A(\blv)$ can be estimated using the trapezoidal rule,
yielding
\begin{equation}\label{eq:A0}
{A}_0
= \frac{\dt}{2}\left(\lambda e^{-0} w_0 + \lambda e^{-\lambda\dt} w_1\right)
= \frac{\dt}{2}\left( w_0 + \lambda e^{-\lambda\dt} w_1\right),
\end{equation}
where $w_0 = \wf[\stv]$ and $w_1=\wf[(\stv-\ev)^+]$,
and $\lambda e^{-\lambda\dt}$ is some constant,
i.e., ${A}_0$ is a weighted sum of $w_0$ and $w_1$.
The update rule for the post-assignment value function $\vf[\stv]$ thus reduces to
\begin{equation}\label{eq:poisson-A-rule}
  \boxed{
    \vf[\stv] \gets
       {A}_0 + e^{-\lambda\dt} \vf[\stv-\ev],
  }
\end{equation}
where $\stv$ has to iterate
the grid points in such order that the updated value for $\vf[\stv-\ev]$
on the right-hand side is always available when updating $\vf[\stv]$.
The lexicographic order
(that can be used for the state-space reduction, see Section~\ref{sec:state-reduction})
is such.
Because $\lambda \dt$ and $e^{-\lambda\dt}$ are constants,
the above expression for $\vf[\stv]$ is simply a weighted sum of
three values, $\wf[\stv]$, $\wf[\stv-\ev]$ and $\vf[\stv-\ev]$.

Note that also Simpson's rule can be applied, e.g., by choosing $h=2\dt$.
Utilizing the lack of memory property in this manner reduces
the number of computational operations significantly as then $\vf[\stv]$
can be obtained in one coordinated ``sweep'' as a weighted sum of
$\wf[\stv]$, $\wf[\stv-\ev]$ and $\vf[\stv-\ev]$
starting from $\stv=0$.
The downside is that
we no longer can take advantage of parallel computations.
(The original update rule is embarrassingly parallel, as each
$\vf[\stv]$ depends solely on $\wf[\stv]$.)

\subsection{Improving value iteration under time-scale separation}
\label{sec:recursive-high-lambda}

As discussed earlier in Section~\ref{sec:time-step}, the time step
$\dt$ is chosen to match the job size distribution $X$
so that, e.g., Simpson's composite rule gives sufficiently accurate results
when expectations over $X$ are computed.
There is the obvious trade-off between accuracy and computational complexity
that while a smaller $\dt$ does improve accuracy,
it also increases the number of grid points and thereby also
the number of computational operations.

\begin{figure}
  \centering
  \sffamily\small
  \begin{tabular}{c@{\hspace{1cm}}c}
  \begin{tikzpicture}[font=\scriptsize,scale=0.65]
    \draw[thin,blue!30] (0,0) grid (9.8,5.2);
    \node at (0,-0.5) {$0$};
    \node at (4,-0.5) {$1$};
    \node at (8,-0.5) {$2$};
    \node at (-0.3, 1) {$1$};
    \node at (-0.3, 2) {$2$};  
    \node at (-0.3, 3) {$3$};
    \node at (-0.3, 4) {$4$};
    \node at (-0.3, 5) {$5$};
    \foreach \la/\col in {1/red,5/blue} {
      \draw[domain=0:9.8, \col!50!black, fill=\col!25, opacity=0.5,smooth, thick, variable=\y]  plot  ({\y}, {\la * exp(-\la*\y/4.0)})
      -- (9.8,0) -- (0,0);
    }
    \node[rotate=-71] at (0.95,2.5) {$A$ with $\lambda=5$};
    \node[rotate=-5] at (4,0.65) {$X$ with $\E{X}=1$};
    \draw[thick] (0,5.2) -- (0,0) -- (9.8,0);
    
  \end{tikzpicture}
  &
  \begin{tikzpicture}[node distance=2mm,scale=0.75]
    \draw[thin,gray!50] (0,0) grid (5,3);
    \draw[thick,<->,>=latex] (0,4.5) --  node [left] {$\bl_2$} (0,0) -- node [below] {$\bl_1$} (6.5,0);
    \draw[line width=4pt,red!60,opacity=0.5] (0,0) -- (2,0) -- (3,1);
    \draw[line width=1pt,black] (3,1) -- (4,2);
    \node[draw,circle,inner sep=2pt,fill,label={above left:$\blv_0$}] at (4,2)  (n0) {};
    \node[draw,circle,inner sep=2pt,fill,label={above left:$\blv_\dt$}]  at (3,1)  (n1) {};
    \node[draw,circle,inner sep=2pt,fill,label={above left:$\blv_{2\dt}$}]  at (2,0) (n2) {};
    %
  \end{tikzpicture}
  \\
  (a) & (b)
  \end{tabular}
  \caption{Left:
    As the system size increases, the time-scales of $A$ and $X$ start to separate,
    which induces numerical challenges.
    Right:
    Updating $\vf(\blv)$ from $\wf(\blv)$ corresponds to taking
    the expectation along the depicted trajectory.
    With Poisson arrivals this can be achieved recursively.
  }
  \label{fig:A-X-scales}
\end{figure}
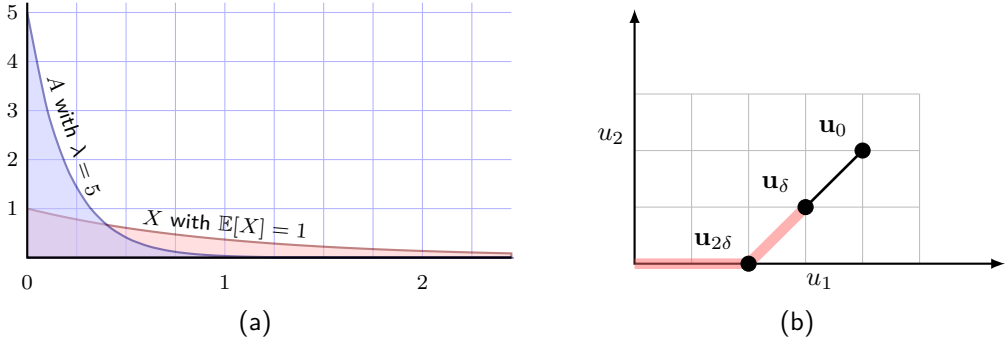

Even when $\dt$ is chosen properly with respect to $X$,
numerical problems are still bound to arise when
the size of the system increases.
Namely, when the number of servers $\sn$ increases
the arrival rate $\lambda = 1/\E{A}$ should also increase to keep the extra servers busy.
In fact, server systems are often scaled so that $\rho$ is kept fixed.
Eventually this leads to a time-scale separation between $X$ and $A$, as
illustrated in Figure~\ref{fig:A-X-scales}.
This means that the time step $\dt$ that was appropriate for $X$
is no longer appropriate when expectations over $A$ are evaluated
(cf.\ Eq.~\eqref{eq:vf-wz2}).
For example, using Simpson's composite rule for the mean gives
$$
\E{A} \approx \frac{\dt}{3} \left( h_0 + 4\,h_1 + 2\,h_2 + \ldots \right), 
\quad\text{where $h_i = i\,\dt\cdot f_A(i\,\dt)$},
$$
and if $f_A(t)$ is very small for $t\ge \dt$,
results tend to become inaccurate. 

\paragraph{Tailored integration:}
Let us consider the integral \eqref{eq:Au} when $h=\dt$.
Instead of applying the trapezoidal rule,
this integral can be computed exactly after
assuming a suitable smooth behavior for $\wf(\blv)$
between the two consecutive grid points $\stv$ and $(\stv-\ev)^+$.
Let $\blv_t$ denote the state when work has been processed from the initial state $\blv$
for time $t$ without new arrivals, i.e., $\blv_t = (\blv-t\ev)^+$.
If we assume that $\wf(\blv_t)$ is linear from $t=0$ to $t=\dt$ along
the given trajectory, then the update rule turns
out to be  an explicit closed-form expression.
Let $w_0=\wf[\stv]$ and $w_1=\wf[\stv-\ev]$. Then
$$
{A}_1 :=
\int_0^{\dt } \lambda  e^{-\lambda  t} \cdot \left(
w_0 + t \frac{w_1-w_0}{\dt}\right) \, dt
=
w_0 - w_1 e^{-\lambda\dt}
+ \frac{\left(1-e^{-\lambda\dt}\right) \left( w_1-w_0 \right)}{\lambda\dt},
$$
yielding, in vector notation,
\begin{equation}
  \label{eq:poisson-A1}
{A}_1 = \frac{1}{\lambda\dt}
\begin{pmatrix}
  w_0 &
  w_1
\end{pmatrix}
\cdot
\left(\begin{array}{rcr}
  \lambda\dt&-&(1-e^{-\lambda\dt}) \\
  -\lambda\dt&+&(1+\lambda\dt)(1-e^{-\lambda\dt})
\end{array}
\right)^T.
\end{equation}
Similarly, the quadratic function $\hat w(t) = at^2+bt+c$
can be fit to values of $\wf(\blv_t)$ at time instants $t=0,\dt,2\dt$,
corresponding to states $\stv$, $\stv-\ev$ and $\stv-2\ev$,
yielding
\begin{equation}
  \label{eq:poisson-A2}
{A}_2 =
\frac{1}{2(\lambda\dt)^2}
\begin{pmatrix}
  w_0 & w_1 & w_2
\end{pmatrix}
\cdot
\left(\begin{array}{rcr}
 - 2 \lambda\dt  (1-\lambda\dt) &+& (2-\lambda\dt) \left(1-e^{-\lambda\dt}\right)  \\
   2  \lambda\dt  (2-\lambda\dt) &-& 2\left(2-(\lambda\dt)^2 \right) \left(1-e^{-\lambda\dt}\right)\\
 - 2 \lambda\dt &+& (2+\lambda\dt) \left(1-e^{-\lambda\dt}\right)
\end{array}
\right)^T.
\end{equation}
Hence, $A_0$, $A_1$ and $A_2$,
given in \eqref{eq:A0}, \eqref{eq:poisson-A1}, and \eqref{eq:poisson-A2}, respectively,
are all weighted sums of $w_0$, $w_1$ and $w_2$ (for $A_2$),
and thus $\vf[\stv]$ can be updated fast in a single coordinated sweep across
all values of $\stv$.

In the numerical examples (see Section~\ref{sec:examples}), 
we observe problems with $\sn=4$ servers when $\rho=0.9$ and $\dt=0.25$.
Using the specifically-tailored integration rule
\eqref{eq:poisson-A2} for one time step
fixes the numerical problem.
Hence, as expected, the specifically tailored integration is
numerically superior to the earlier versions
when $\lambda$ increases together with the system size $\sn$.
We note that 
\eqref{eq:poisson-A1} and
a smaller $\dt$ would also help to tackle the numerical problems.
However,
we chose to use \eqref{eq:poisson-A2} as
the quadratic form because
it is more accurate than the linear relation used in \eqref{eq:poisson-A1} with
similar computational complexity.
Smaller $\dt$, on the other hand,
would either lead to a larger state space or earlier truncation of the backlogs.

\subsection{State-space reduction}
\label{sec:state-reduction}

\begin{figure}
  \centering
  \includegraphics[width=14cm]{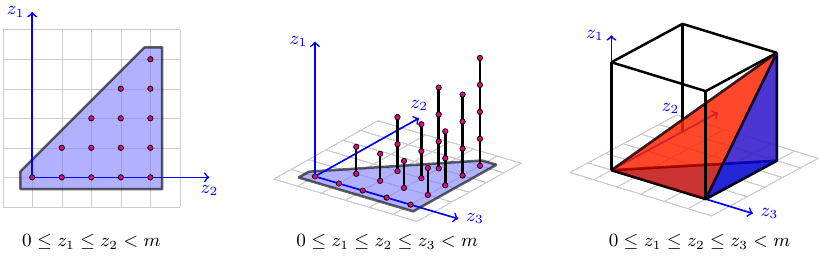}
  \caption{State-space reduction with $m=2$ and $m=3$ servers.}
  \label{fig:state-space-3d}
\end{figure}

Discretizing the truncated state space yields a grid with $\gs^\sn$ points.
However, due to symmetry, it is sufficient to consider only points
$\stv = (\st_1,\ldots,\st_{\sn})$ with $0 \le \st_{1} \le \ldots  \le \st_{\sn} <\gs$.
Figure~\ref{fig:state-space-3d} illustrates this for $\sn=2$ and $\sn=3$ servers.

The size of the original state space is $\gs^\sn$, and
the size of the reduced state-space (i.e., the number of unknown variables
after the reduction)
is given by \cite{hyytia-peva-2024}
\begin{equation}\label{eq:ordered-size}
  ||\mathcal{S}|| = \frac{ \gs(\gs+1)\cdots (\gs+\sn-1)}{ \sn!}.
\end{equation}
The grid points can be enumerated using the lexicographic order
so that the rank of grid point $\stv$ is
$$
\text{pos}(\stv) = \binom{\st_1}{1} + \binom{\st_2+1}{2} + \ldots + \binom{\st_{\sn}+\sn-1}{\sn},
$$
which serves as a linear index to the memory array where the multi-dimensional grid is stored.

\section{Numerical experiments}
\label{sec:examples}

Let us next take a look at some numerical results for size-aware dispatching systems with $\sn=2,\ldots,6$ servers,
Poisson arrival process and exponentially distributed job sizes.
These results complement the numerical results given in \cite{hyytia-peva-2024}.

\subsection{Optimal policy for two servers}

\begin{figure*}
  \centering
  \scriptsize
  \begin{tabular}{lcccc}
    & Value function $\vf(\bl_1,\bl_2)$ & Contour lines & Diagonal and Boundary & Fixed $\bl_2$ \\
    \rotatebox{90}{\hspace{15mm}$\rho=0.4$} &
    {\includegraphics[width=42mm,trim=0 0 0 40]{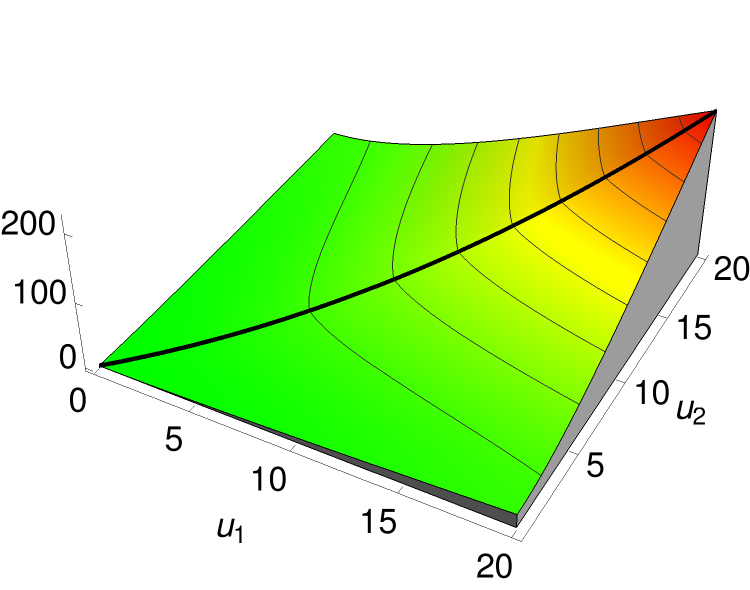}}   &
    \includegraphics[width=35mm]{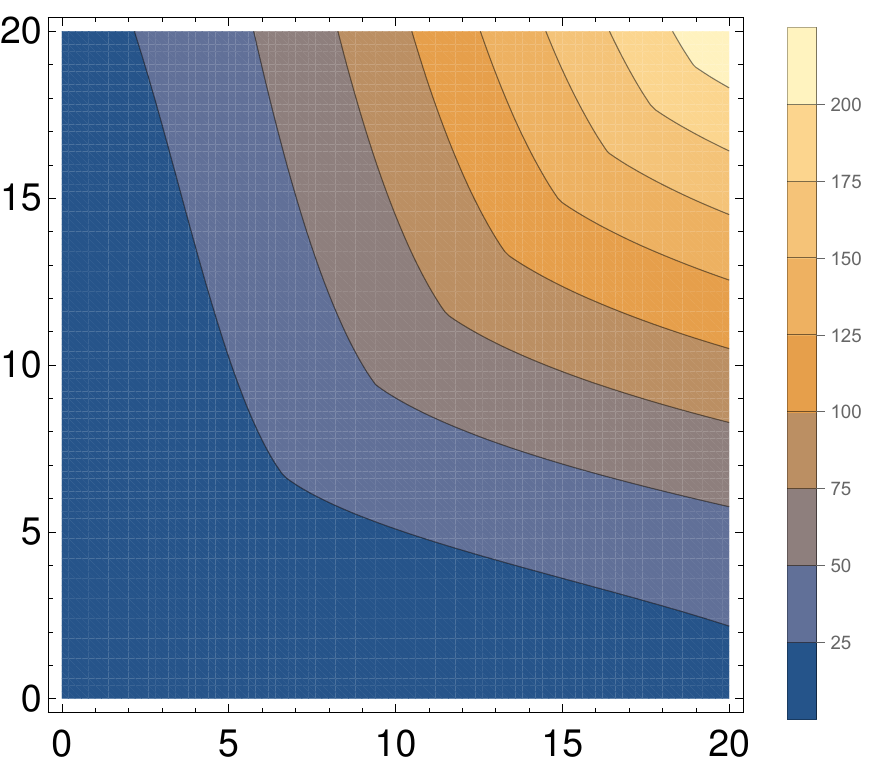} &
    \includegraphics[width=32mm]{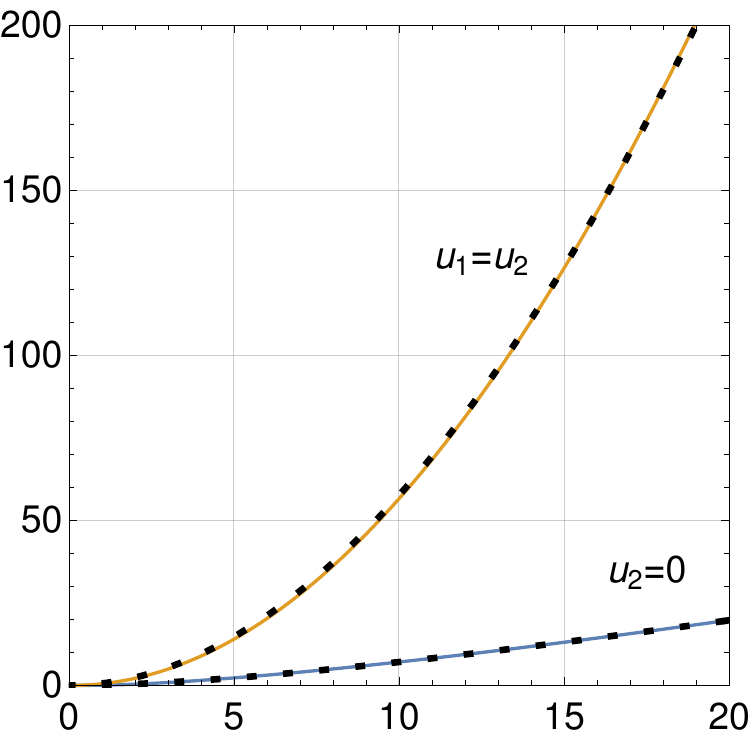} &
    \includegraphics[width=32mm]{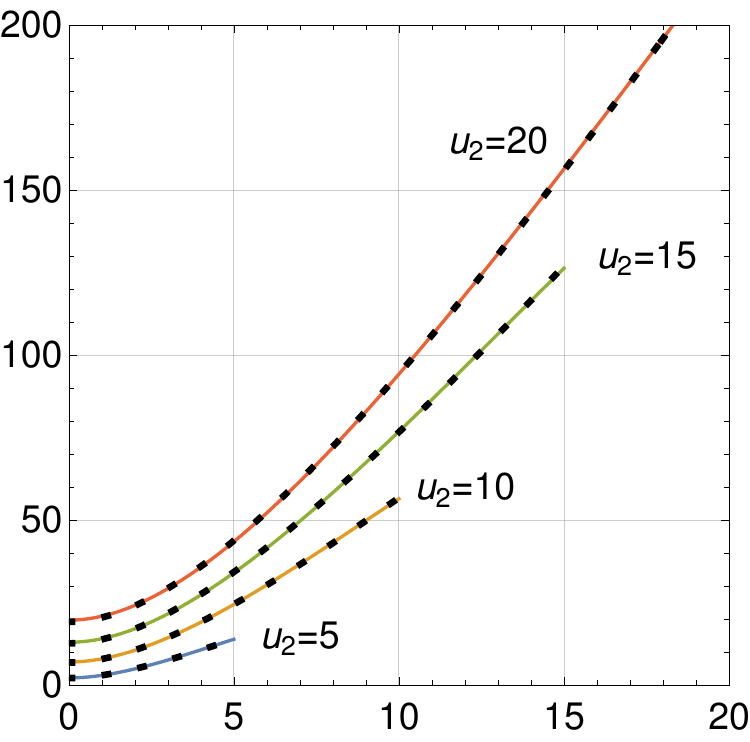} 
    \\
    \rotatebox{90}{\hspace{15mm}$\rho=0.9$} &
    {\includegraphics[width=42mm,trim=0 0 0 40]{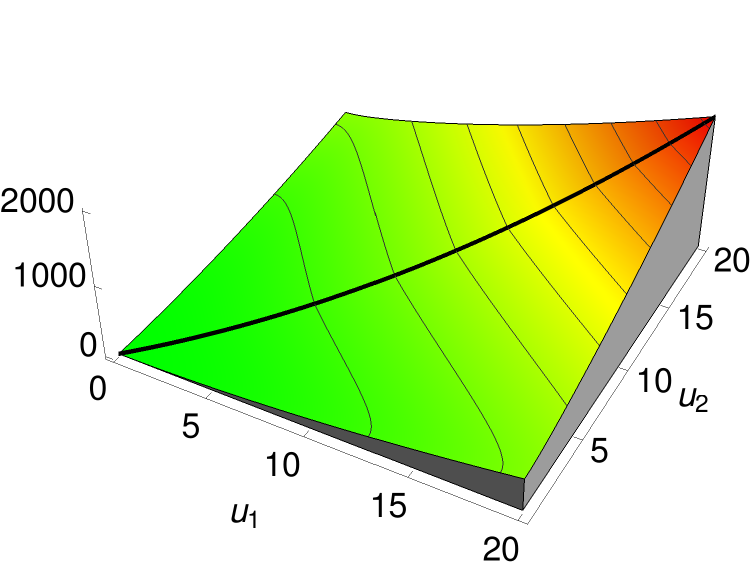}} &
    \includegraphics[width=35mm]{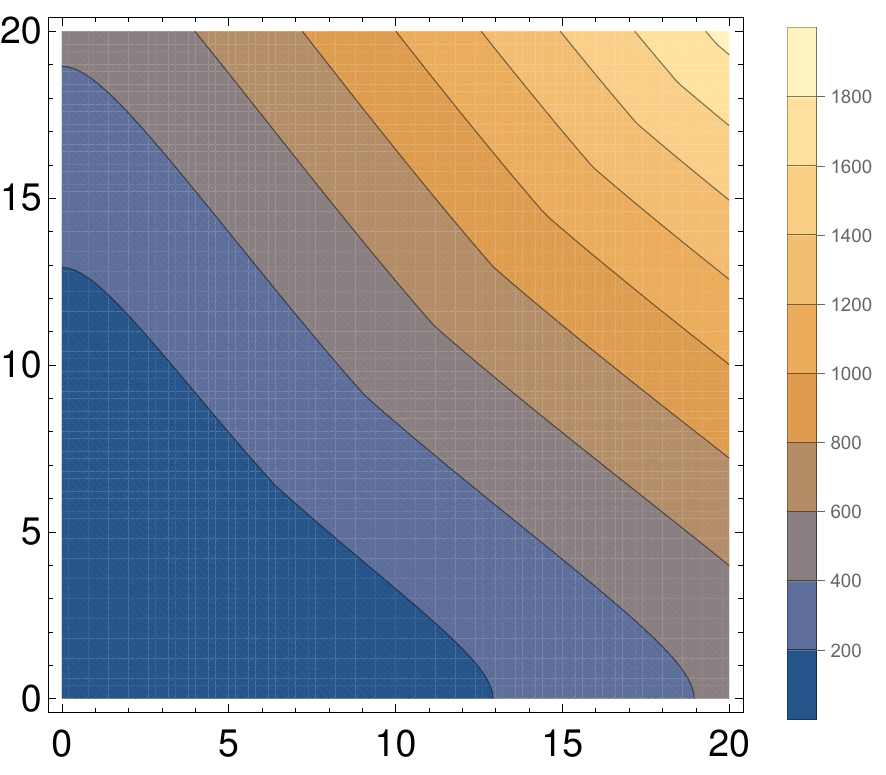} &
    \includegraphics[width=32mm]{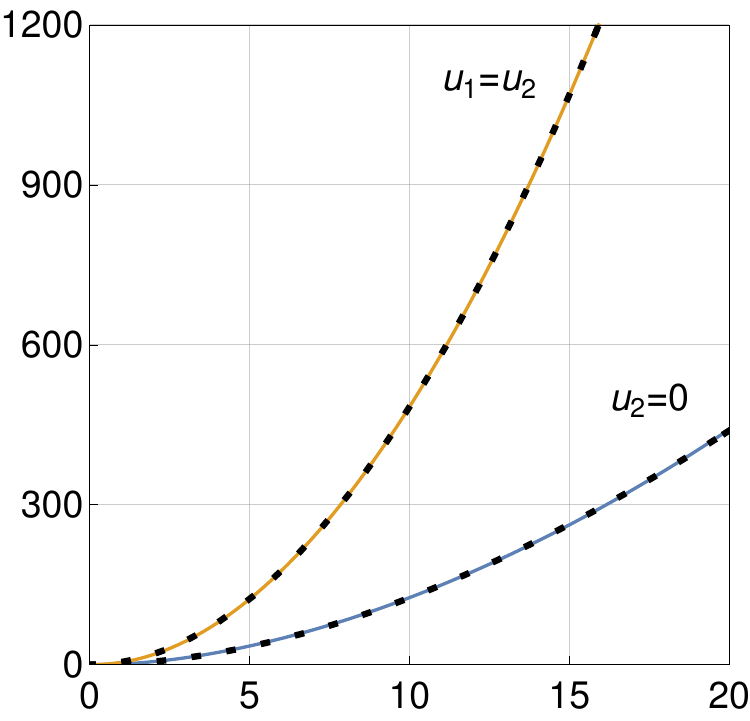} &
    \includegraphics[width=32mm]{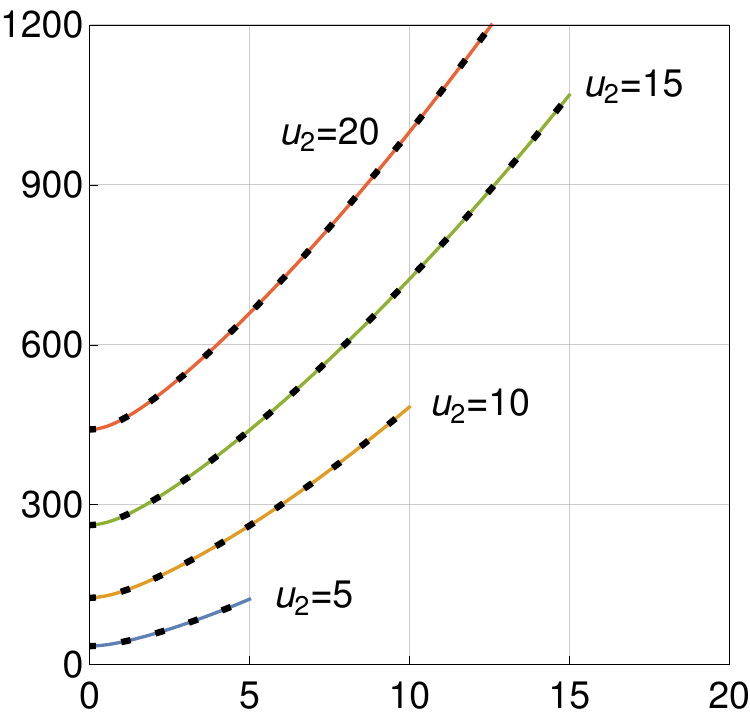} 
  \end{tabular}
  \caption{The value function for two servers with $\rho=0.4$ and $\rho=0.9$
      when job sizes are exponentially distributed.
      The black dots correspond to functions fitted to match the numerical values.
      For $\rho=0.4$,
      the value function along the diagonal ($\bl_1=\bl_2$) appears to be nearly quadratic,
      ${\vf(\bl,\bl)\approx 0.56 u^2}$.
      For all other cuts, we found good approximations
      of form $a + b u^2 /(c+u)$, where $a$, $b$, and $c$
      are properly chosen.}
    \label{fig:contour}
\end{figure*}

The interesting question is how these value functions actually look, i.e., what kind of shape
do they have?
Figure~\ref{fig:contour} depicts the value function for the two-server system
with Poisson arrivals and exponentially distributed job sizes when $\rho=0.4$ and $\rho=0.9$,
corresponding to relatively low and high loads.
As before, the objective is to minimize the mean waiting time.
The big picture is that (these two) value functions change smoothly and are approximately quadratic
along the diagonal,
and easy to approximate along the shown other cuts
using elementary functions of form $a + bu^2/(c+u)$,
as shown in the figures on the right (the black dots correspond to the fitted functions and
solid curves to the numerical results).
However, the exact shape must be determined
by solving an optimization problem.
When $\rho=0.4$ the equivalue contour lines curve strongly and seem to align with the axes as the backlogs increase.
However, as the load increases the contour lines seem to approach straight lines of form
$\bl_1 + \bl_2 = c$, where $c$ is some constant.

\begin{figure*}
  \centering
  \scriptsize
  \begin{tabular}{c@{\hspace{20mm}}c}
    $\rho=0.4$ & $\rho=0.9$\\
    \includegraphics[width=34mm]{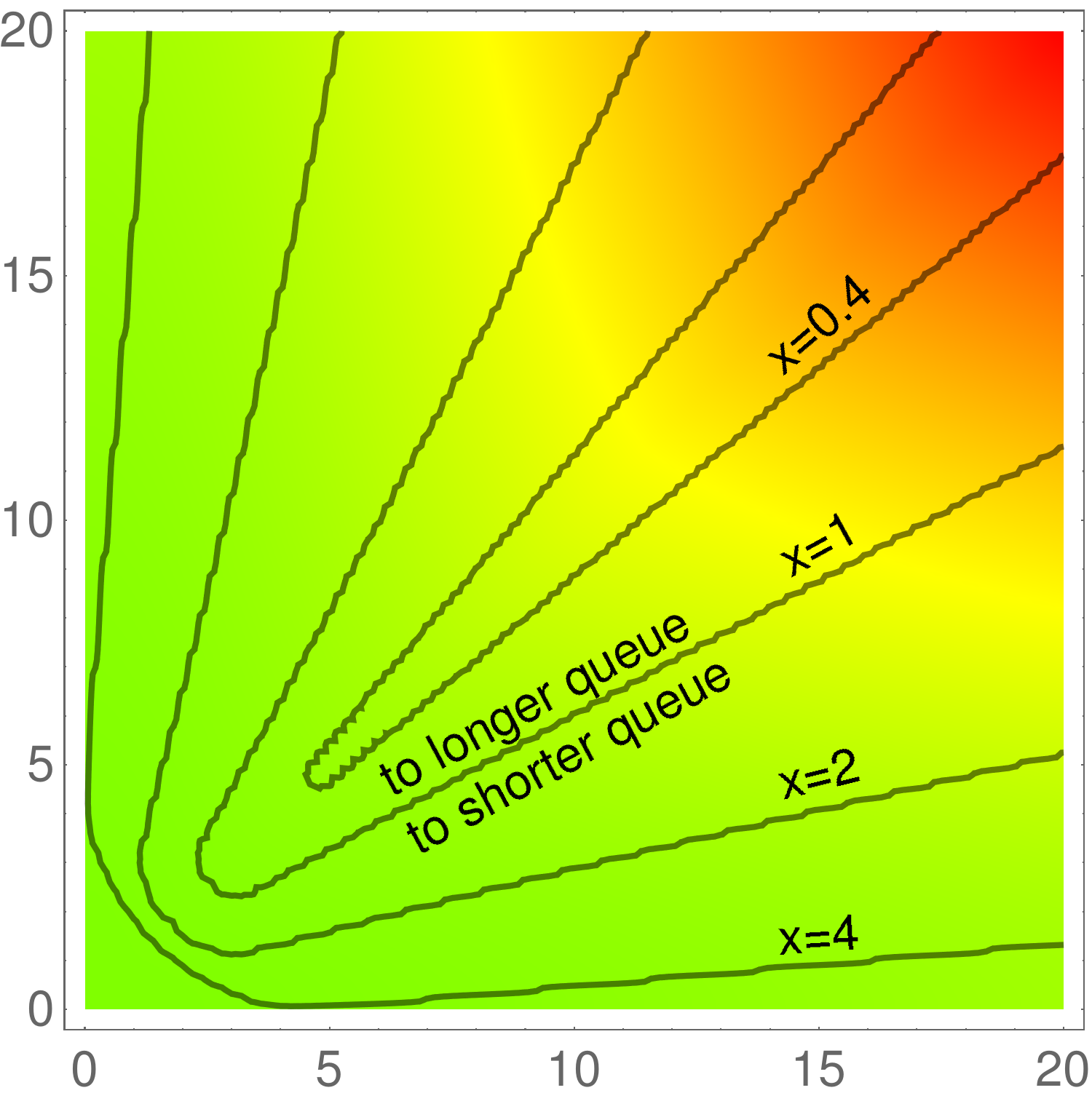} &
    \includegraphics[width=34mm]{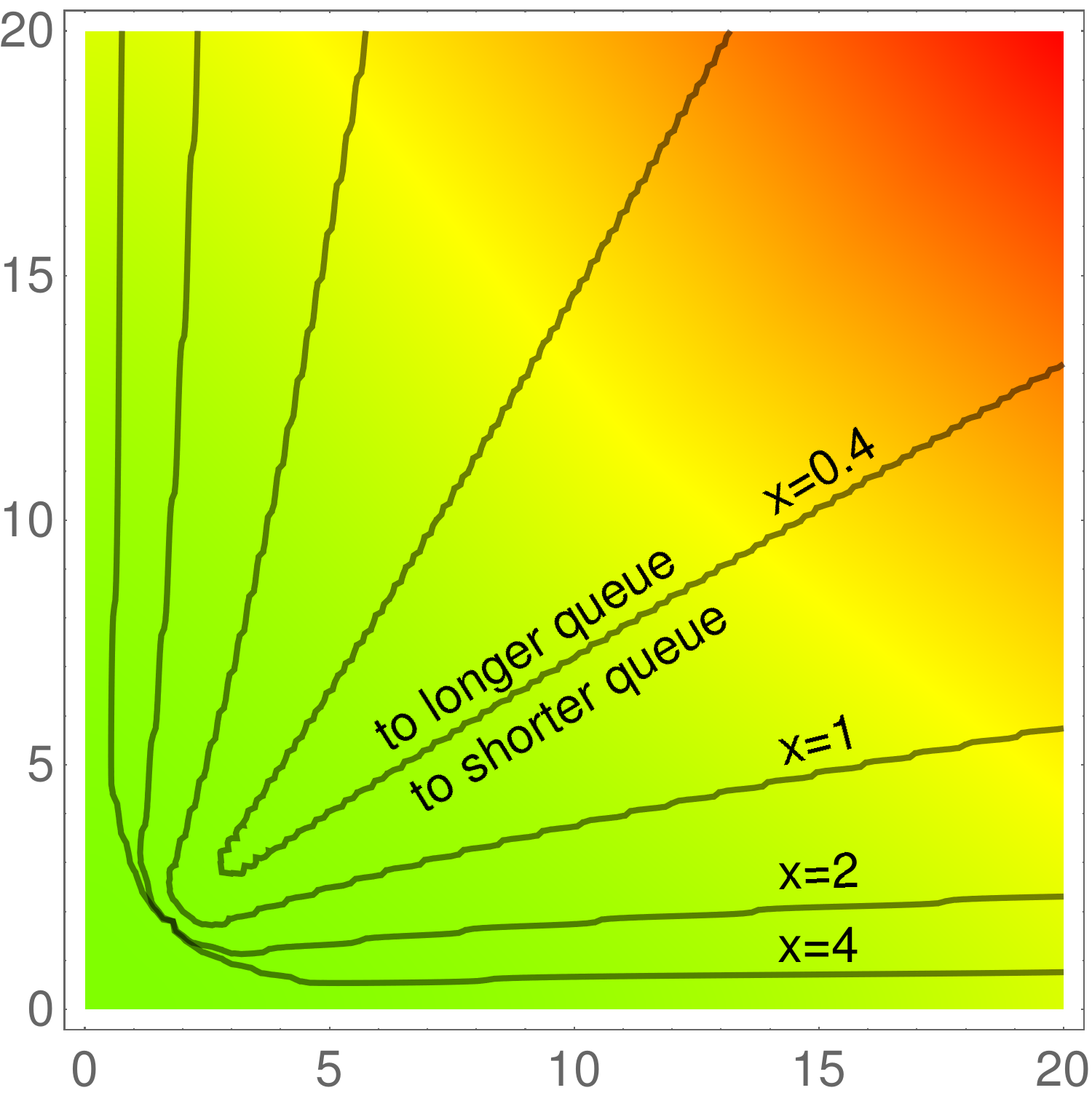}
  \end{tabular}
  \caption{Optimal actions with two identical servers and exponentially
      distributed jobs for size $x$ job, $x \in \{0.4,\,1,\,2,\,4\}$.
      In the interior region, the job is forwarded to the longer queue,
      and vice versa.} 
  \label{fig:opt-action-exp2}
\end{figure*}

\begin{figure}
  \centering
  \begin{tabular}{cccc}\scriptsize
  \parbox{34mm}{
    \begin{tikzpicture}[scale=0.48,font=\tiny]
      \foreach \i in {1,...,5} {
        \draw[magenta!60!black] (0,\i) -- node [above,sloped,black] {$u=\i$} (\i,0);
      }
      \draw[<->,>=latex,thick] (0,6) node [left,sloped] {$u_2$} -- (0,0) -- (6,0) node [below] {$u_1$};
    \end{tikzpicture}
    \\
    \centerline{\scriptsize Diagonal trajectories}
  }
  & \rotatebox{90}{\hspace{-3mm}$\rho=0.4$} 
  & \parbox{65mm}{\includegraphics[width=56mm]{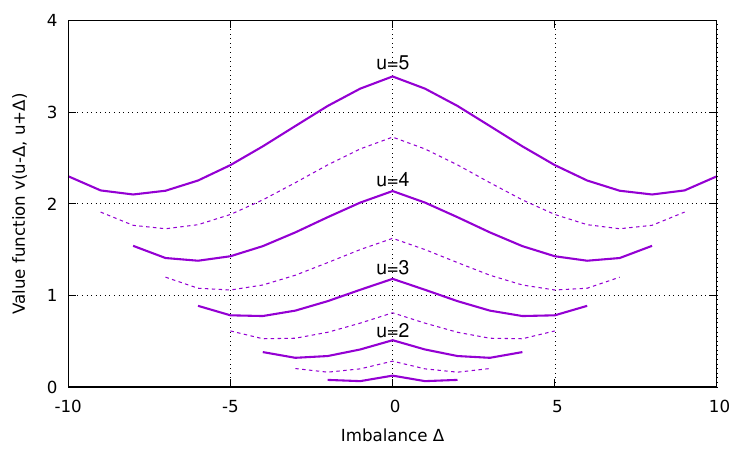}}
  & \parbox{65mm}{\includegraphics[width=56mm]{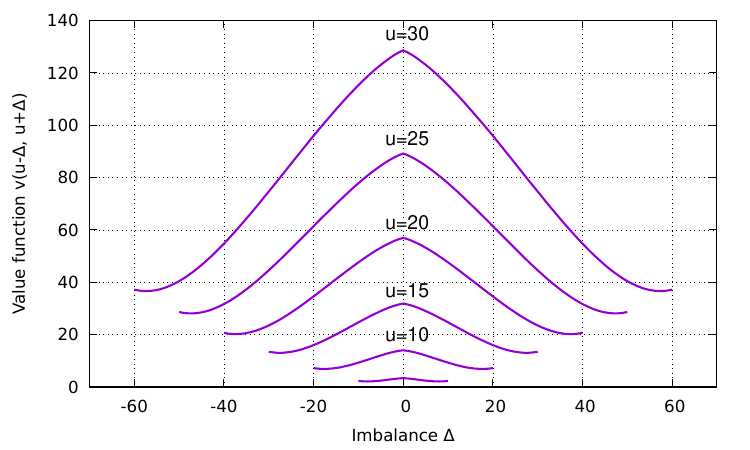}}\\
  & \rotatebox{90}{\hspace{-3mm}$\rho=0.9$} 
  & \parbox{65mm}{\includegraphics[width=56mm]{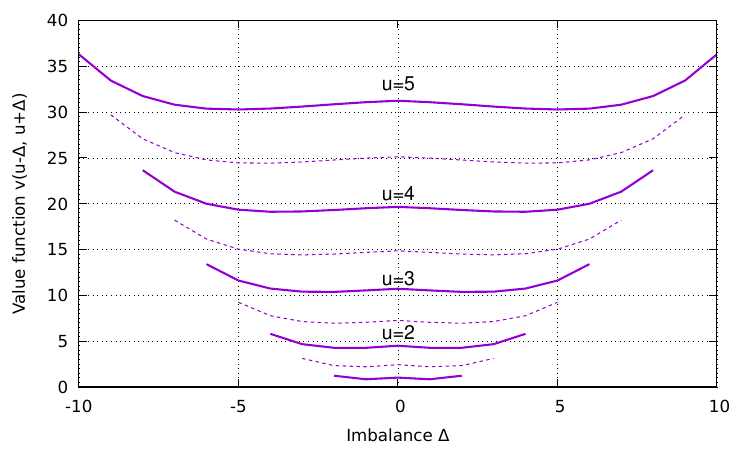}}
  & \parbox{65mm}{\includegraphics[width=56mm]{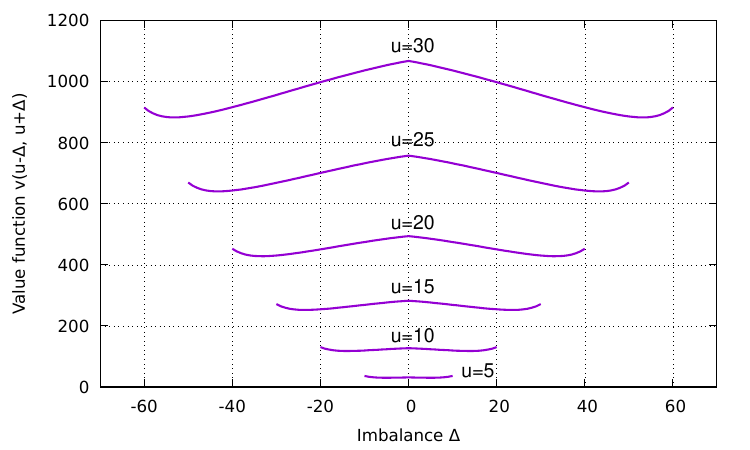}}
  \end{tabular}
  \caption{Value function near the origin (middle) and in larger scale (right)
    depicted for fixed total backlogs $u=u_1+u_2$ for $\rho=0.4$ (upper) and $\rho=0.9$ (lower).
    The fixed total backlogs correspond to the
    ``diagonal trajectories'' as shown on the top left.}
  \label{fig:diag}
\end{figure}
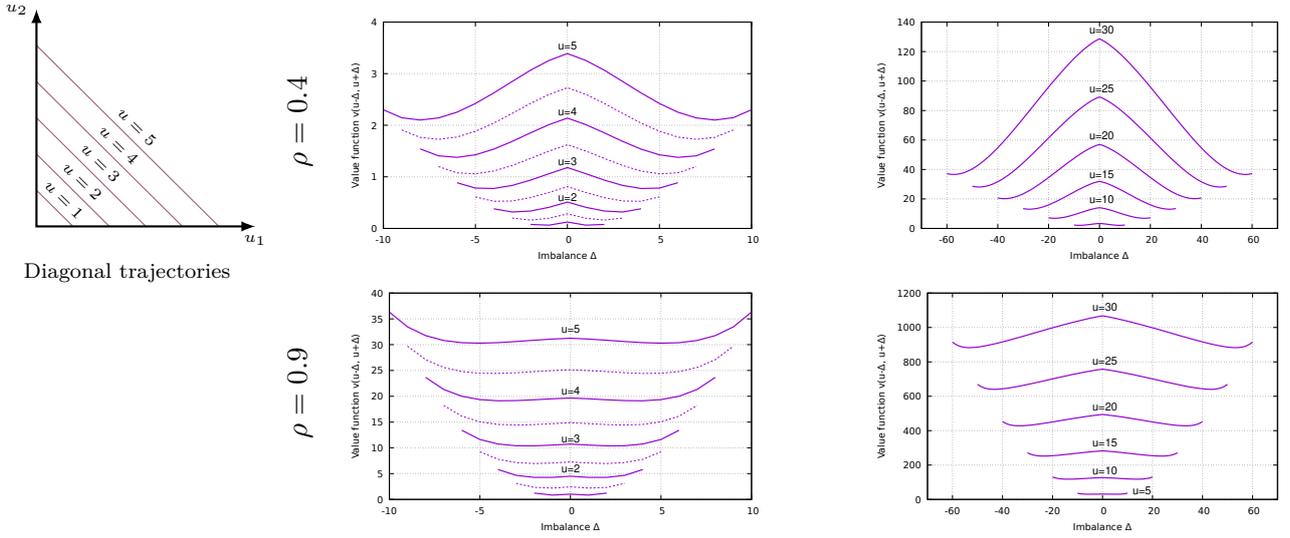

Consequently, let us next take a closer look at the same two value functions along the cuts where
the total backlog $u=u_1+u_2$ is fixed and some multiple of the discretization step $\dt$ (here $\dt=0.25$).
Then let $\Delta = i\dt$ denote the imbalance,
$\Delta = \bl_2 - \bl_1$. 
Value functions across such cuts go via grid points and are symmetric with respect to $\Delta$,
yielding different ``moustache'' patterns shown in Figure~\ref{fig:diag}.

Figure~\ref{fig:opt-action-exp2} depicts the optimal dispatching actions
in the two cases for jobs with size $x \in \{0.4,\, 1,\, 2,\, 4\}$.
It turns out that the optimal policy is characterized by size-specific
curves that divide the state-space into two regions: in the outer region (which includes
both axes) jobs are routed to the shorter server (to avoid idling), and in the interior
region actions are the opposite (to keep the servers unbalanced).
  Note that, e.g., short jobs are
  assigned to the shorter queue almost always, and vice versa.

\subsection{Optimal policy for three servers}

Let us next consider a slightly larger system of three servers.
In this case, we assume that jobs arrive according to Poisson process with rate $\lambda=2.7$,
job sizes are exponentially distributed, $X \sim \text{Exp}(1)$,
and the offered load is thus high, $\rho=0.9$.

\begin{figure}
  \centering
  \scriptsize
  \begin{tabular}{cccc}
    $x=0.5$ &
    $x=1$ &
    $x=2$ &
    $x=3$ \\
    \includegraphics[width=3.5cm]{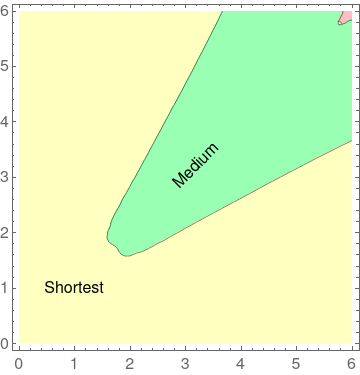} &
    \includegraphics[width=3.5cm]{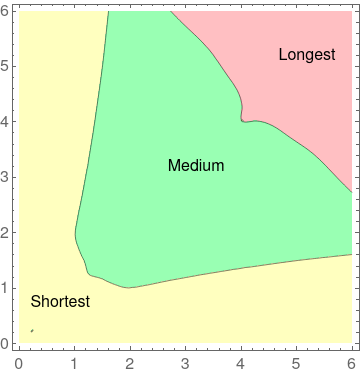}  &
    \includegraphics[width=3.5cm]{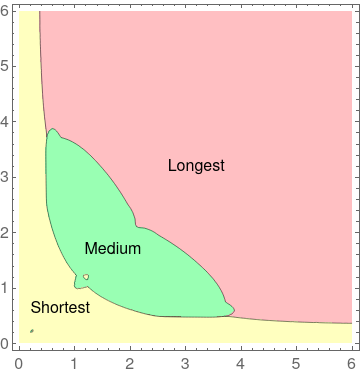}  &
    \includegraphics[width=3.5cm]{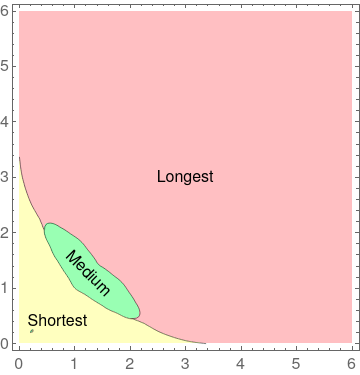}
  \end{tabular}
  \caption{Optimal dispatching actions with $\sn=3$ servers for size $x$ job when the backlog in the longest server is
    fixed to $\bl_3=8$
    and the offered load is $\rho=0.9$.
    Figures may have some artifacts due to the `discretization'' when $\bl_1=\bl_2$ and the two shortest queues are
    equivalent.}
  \label{fig:aktion-3}
\end{figure}

Illustrating the actions of optimal policy becomes non-trivial as soon as the number of servers increases.
Figure~\ref{fig:aktion-3} shows the optimal actions as a function of
the shorter backlogs
$\bl_1$ and $\bl_2$
when the backlog in the longest queue $3$ is fixed to $\bl_3=8$.
The size of the new job varies from $x=0.5$ (left) to $x=3$ (right).
We can see the expected pattern
that shorter jobs are again routed to shorter queues
and longer jobs to longer queues.
Boundaries between the different actions may contain small errors due to discretization and 
our visualization method. In particular, when $\bl_1 = \bl_2$, the two servers are equivalent.

\subsection{Convergence of the value iteration}

One important question about the value iteration is how fast it converges,
i.e., how many iteration rounds are needed to obtain a (near) optimal policy.
In our setting, after the state-space truncation,
the number of variables grows exponentially as a function
of the number of servers $\sn$, as given by \eqref{eq:ordered-size}.
The minimum operator over the servers in \eqref{eq:dynprog-ggk3} 
means that our expressions yield a system of non-linear equations.
Let us again consider the familiar example with Poisson arrival process
and $\text{Exp}(1)$-distributed job sizes.
The offered load is $\rho=0.9$ unless otherwise stated.
We shall compare two variants of value iteration.
The first uses the standard Simpson's method when evaluating $\vf[\stv]$.
The second variant utilizes our tailored approximation based on
the lack of memory property of the Poisson process and Eq.~\eqref{eq:poisson-A2}
in particular.
We refer to these two as {\em basic} and {\em w2}.
Note that results with \eqref{eq:poisson-A1} are very similar.

We study two metrics that capture the convergence of the value iteration process.
The first metric is the mean squared adjustment (or change) defined as
\begin{equation}\label{eq:metric-E}
E := \frac{1}{N} \sum_{\stv} \left( \vf_{j+1}[\stv] - \vf_{j}[\stv] \right)^2,
\end{equation}
where $N$ denotes the number of grid points.
Given the iteration converges, $E \to 0$.
The second metric is the estimated mean waiting time $\E{W}=w_0$,
\begin{equation}\label{eq:metric-w0}
w_0 = \E{ \vf( X \ev_1 ) },
\end{equation}
that we obtain automatically during the iteration.
Value iterations start from the random split (RND). 

\paragraph{Two servers:}
Let us start with $\sn=2$ servers.
Figure on the left depicts $E$ and the figure on the right
the estimated mean waiting time $\E{W}$.
Note that the $x$-axess are in logarithmic scale.

\begin{center}
  \includegraphics[width=72mm]{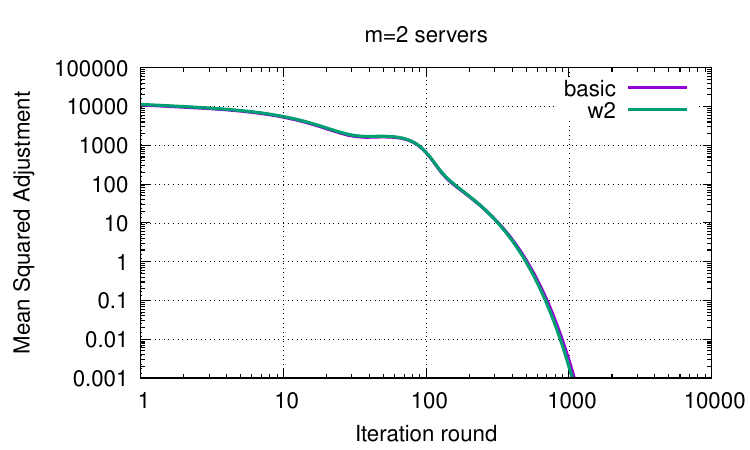}
  \quad
  \includegraphics[width=72mm]{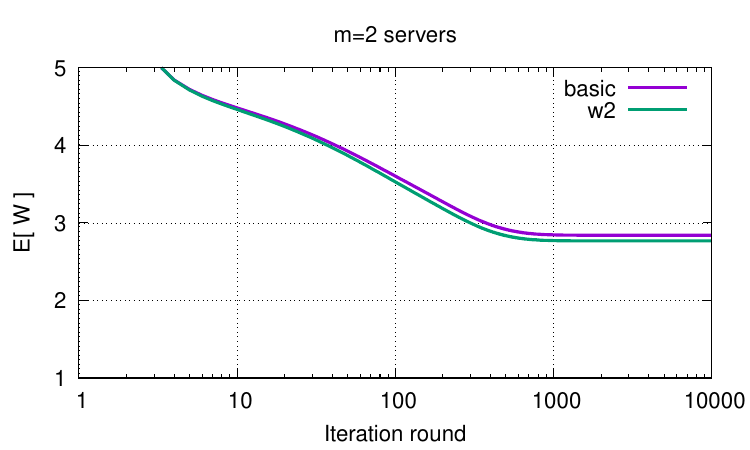}
\end{center}

For $\E{W}$ we can identify three phases: initially during the first 3-4 rounds
$\E{W}$ decreases rapidly. After that $\E{W}$ decreases slowly
for about 1000 rounds, after which the iteration has converged.

For $E$ we can the situation is similar. During
the first 100 rounds values change significantly as
the shape of the value function is not there yet.
After than $E$ decreases rapidly towards zero.
We can conclude that the iteration has practically
converged at this point, say after 1000 rounds.

\paragraph{Three servers:}
The figure below shows the convergence in the equivalent high load setting with $\rho=0.9$ for $\sn=3$ servers.
The numerical results are very similar, but we now notice that the basic algorithm already gives a different
estimate for the mean waiting time $\E{W}=w_0$.
Note also that $\E{W}$ is smaller than with $\sn=2$ servers, as expected.
\begin{center}
  \includegraphics[width=72mm]{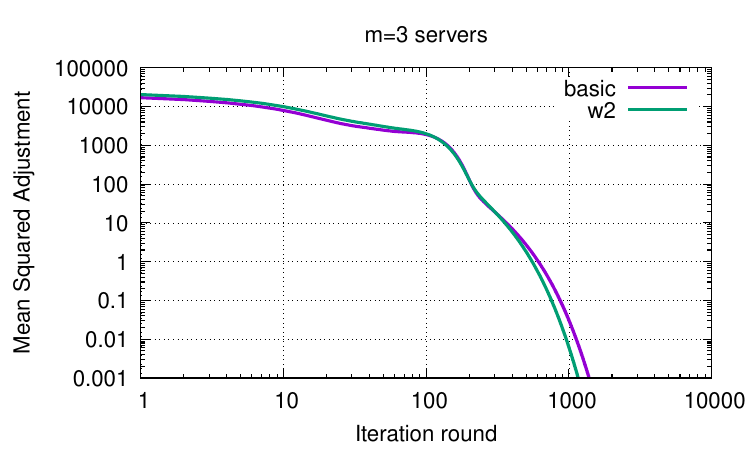}
  \quad
  \includegraphics[width=72mm]{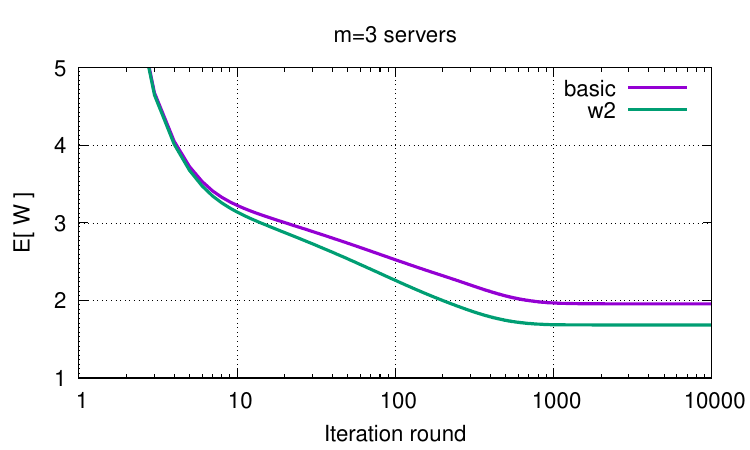}
\end{center}

\clearpage
\paragraph{Four servers:}

With four servers the basic approach using Simpson's composite rule practically fails as the estimated
mean waiting time $\E{W}$ starts to increase after 100 iteration rounds!
In contrast, the tailored approach (w2) based on explicit integration for the one time step when
computing $\vf[\stv]$ from $\wf[\stv]$ still works well.
With the latter, both curves share the same overall pattern as we observed
with $\sn=2$ and $\sn=3$ servers.
\begin{center}
  \includegraphics[width=72mm]{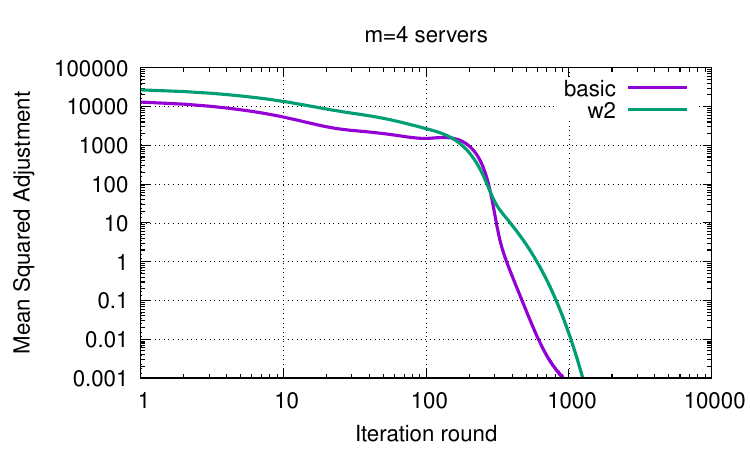}
  \quad
  \includegraphics[width=72mm]{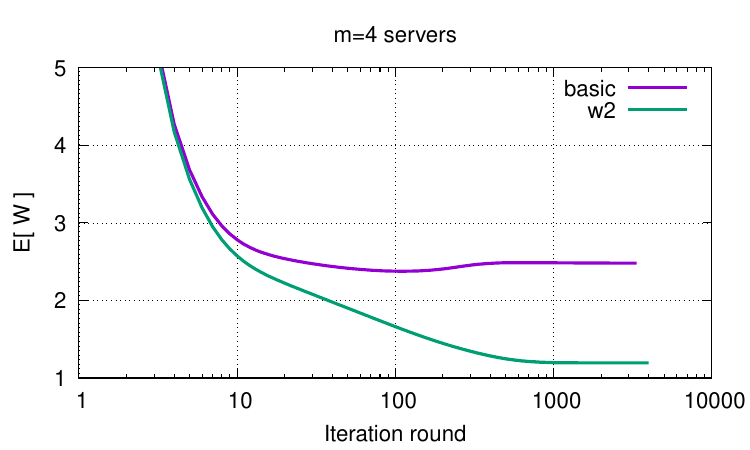}
\end{center}

\paragraph{Comparison of two, three and four server systems:}
The figures below depict the situation with $\sn\in\{2,3,4\}$ servers.
The solid curves correspond to the basic algorithm, and
dashed curves use our tailored recursion (w2).
We can see that in all cases the value iteration converges in about 1000 rounds,
i.e., the number of servers seems to have little impact on this quantity, at least
when the number of servers is relatively small.
The numerical problems that the basic algorithm encounters with $\sn=4$ servers are also clearly visible
in the figure on the right.
\begin{center}
  \includegraphics[width=72mm]{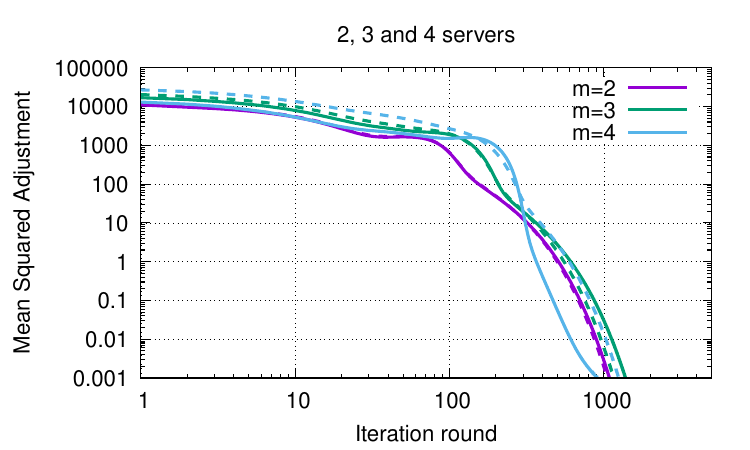}
  \quad
  \includegraphics[width=72mm]{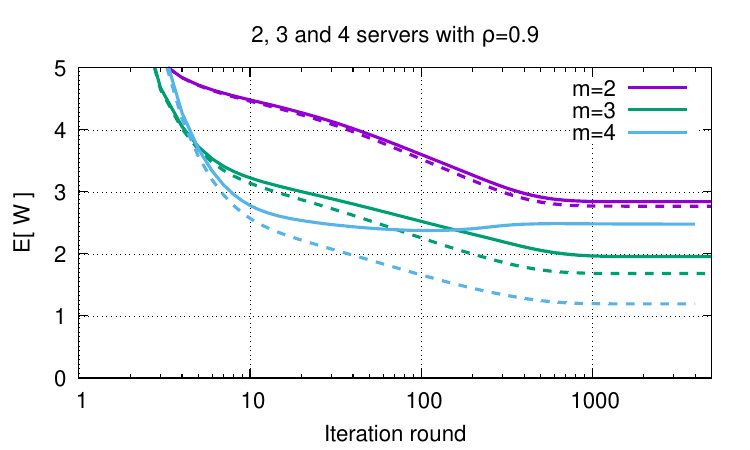}
\end{center}

\paragraph{Larger systems:}

Now we consider a bit larger systems with $\sn=5$ and $\sn=6$ servers.
In this case, we have to decrease the grid size parameter $\gs$ from $200$ to
$120$ ($\sn=5$) and $70$ ($\sn=6$), respectively.
With these values, the size of the value function data structure,
after state space reduction, is slightly below $2$\,GB.
The convergence plots shown below are similar to earlier ones.
Due to a smaller $\gs$,
the convergence in terms of iteration rounds is actually slightly faster
(even though the wall-clock running time is significantly longer).
\begin{center}
  \includegraphics[width=72mm]{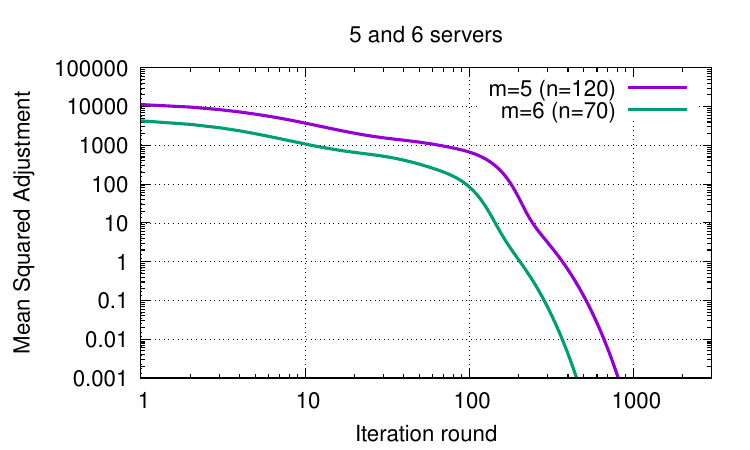}
  \quad
  \includegraphics[width=72mm]{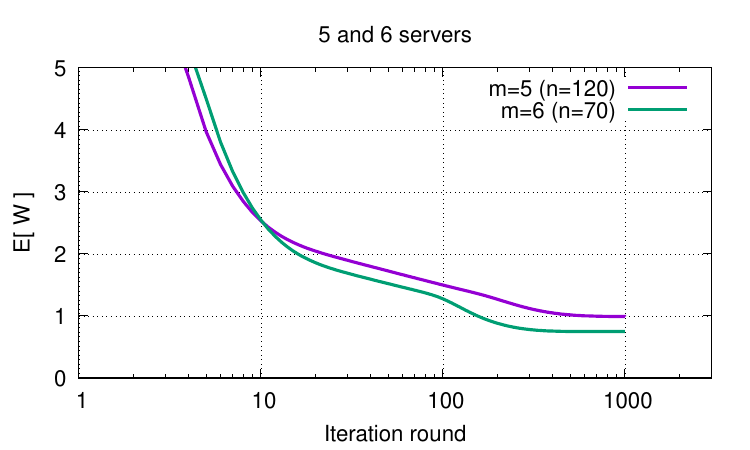}
\end{center}

\clearpage

\subsection{Starting point for value iteration}

Next we study how the initial value function $\vf_0(\blv)$
affects the convergence. We consider three options:
\begin{enumerate}
\item {\bf From zero:} The most simple option is to set $\vf_0(\blv)=0$ for all $\blv$. 
\item {\bf From RND:} Dispatching systems with Poisson arrival process and random split
  decompose into $\sn$ statistically identical M/G/1 queues.
  Given the value function for the M/G/1 queue is available in closed form \cite{hyytia-ejor-2012,hyytia-peva-2020},
  the decomposition gives us the value function for the whole system,
  $$
  \vf_0(\blv) = \sum_{i=1}^{\sn} \frac{\lambda' (\bl_i)^2}{2(1-\rho)},
  $$
  where $\lambda'$ is the server-specific arrival rate, $\lambda'=\lambda/\sn$.
\item {\bf From $\boldsymbol{\rho}'$:} We can also start from a value function that has been obtained (numerically) for a different load $\rho'$
  and set $\vf_0(\blv) = \vf^{(\rho')}(\blv)$ for all $\blv$.
\end{enumerate}

\begin{figure}
  \centering
  \includegraphics[width=7cm]{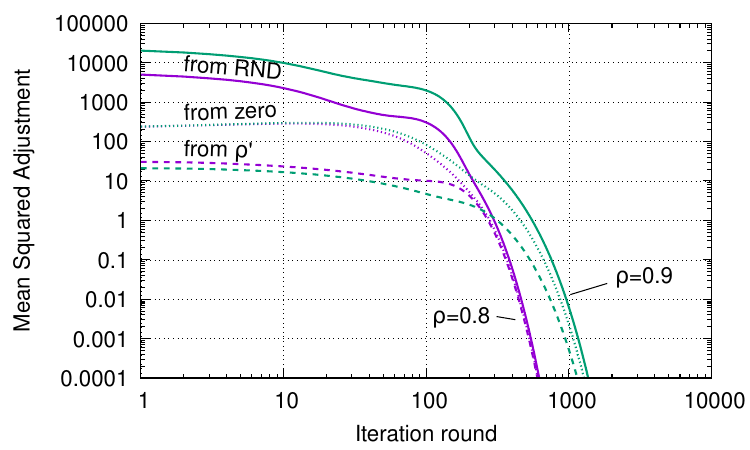}
  \hspace{1cm}
  \includegraphics[width=7cm]{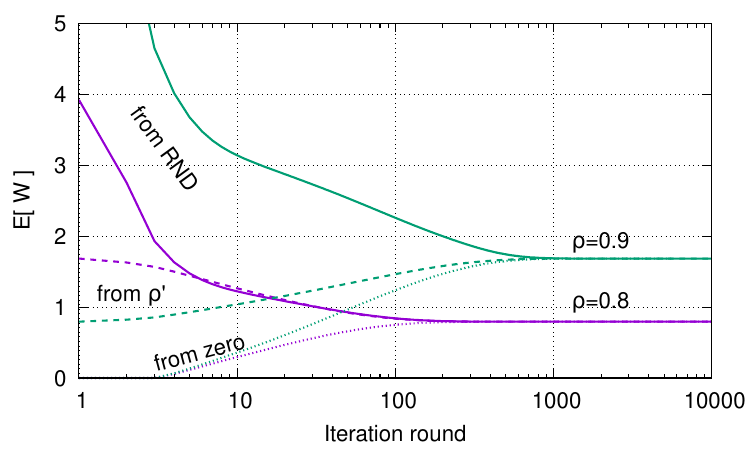}
  \caption{Convergence of the value iteration for $\sn=3$ servers with different initial values
    i)   $\vf_{0}[\stv] =  \vf^{(\text{RND})}[\stv]$ (value function of RND policy),
    ii)  $\vf_{0}[\stv] =  0$ (all values initialized to zero), and
    iii) $\vf_{0}[\stv] =  \vf^{(\rho')}[\stv]$ (value function determined for another $\rho$).
    The difference in convergence time is negligible.}
  \label{fig:vf-inits}
\end{figure}

Here we consider two cases: $\rho_1=0.8$ and $\rho_2 = 0.9$,
and determine the optimal value function for both starting from the different initial values.
For the $\rho'$-method we use the value function of $\rho=0.8$ when determining the value function
for $\rho=0.9$, and vice versa.
The numerical results are depicted in Figure~\ref{fig:vf-inits}.
Somewhat surprisingly, the initial values for value iteration seem to have little impact
to the convergence. Even the optimal value function for a slightly different $\rho$
does not seem to speed-up computations (even though the initial adjustments are smaller).
In this sense, given we are prepared to iterate until the values converge,
one may as well initialize the multi-dimensional array with zeroes.

\subsection{Fairness and the optimal dispatching policy}

The optimal dispatching policies we have determined minimize the mean waiting time.
This is achieved by giving shorter jobs dedicated access to shorter queues.
In particular, routing decisions ensure that backlog is short at least in one queue.
That is, the backlogs are kept unbalanced and only short jobs are allowed to
join the shorter queues.

Such actions tend to introduce unfairness when jobs are treated differently based on their sizes.
To explore this numerically, we have simulated the system with two and three servers and
collected statistics on (waiting times, job size)-pairs.
This data enables us to determine empirical conditional waiting time distributions.
We also simulated the two systems with random split (RND) and the least-work-left policy (LWL).

\begin{figure}
  \centering
  \begin{tabular}{cc}
    \scriptsize Two servers &
    \scriptsize Three servers\\
  \includegraphics[width=76mm]{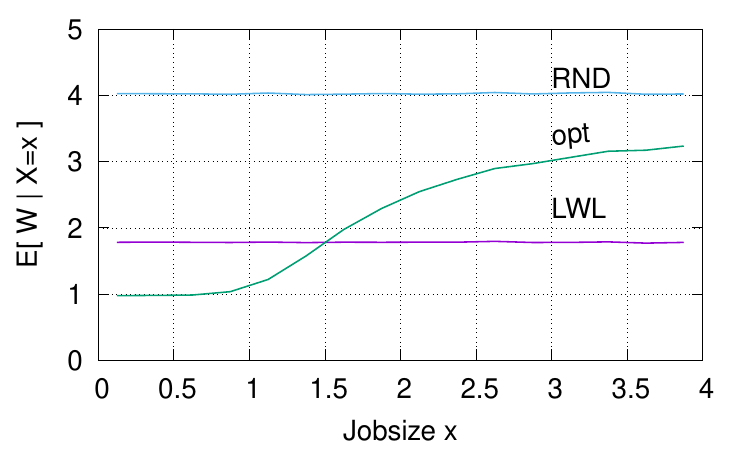} &
  \includegraphics[width=76mm]{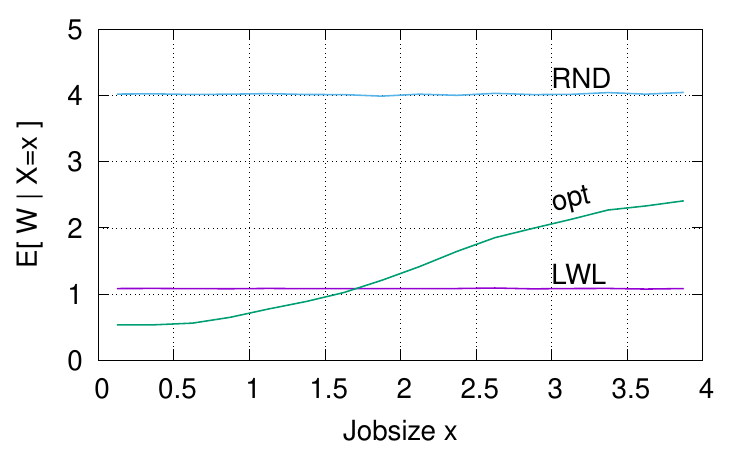}
  \end{tabular}
  \caption{Mean waiting time conditioned on the job size $X$ with RND, LWL and the (near) optimal policy.
      (two and three servers, $X\sim\text{Exp}$, and the offered load is $\rho=0.8$).}
  \label{fig:cond-wt}
\end{figure}

Figure~\ref{fig:cond-wt} depicts the mean waiting time conditioned on the job size
for two- and three-server systems with exponentially distributed job sizes.
As we know, RND and LWL are oblivious to the size of the current job. 
In contrast, the optimal policy clearly gives higher priority for shorter jobs in the same fashion
as under the SRPT and SPT scheduling disciplines.
This type of service differentiation is easier if the queues are unbalanced, which is
the reason why the optimal dispatching policy actively
seeks to maintain unequal queue lengths.

\begin{figure}
  \centering
  \begin{tabular}{ccc}
    \scriptsize $\rho=0.6$ &
    \scriptsize $\rho=0.9$ \\
  \includegraphics[width=76mm]{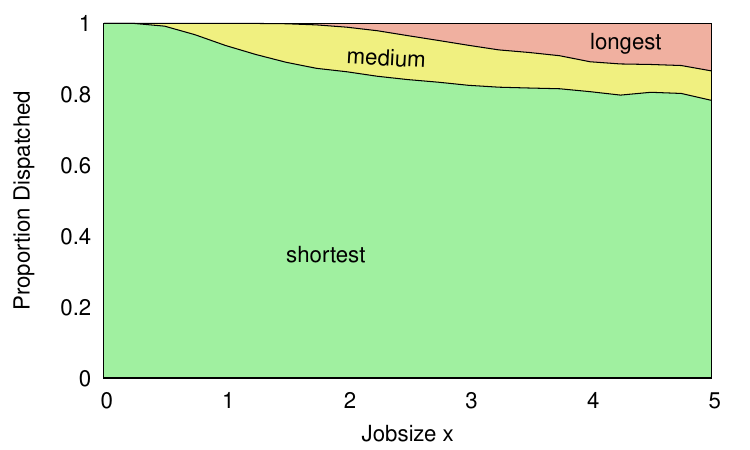} &
  \includegraphics[width=76mm]{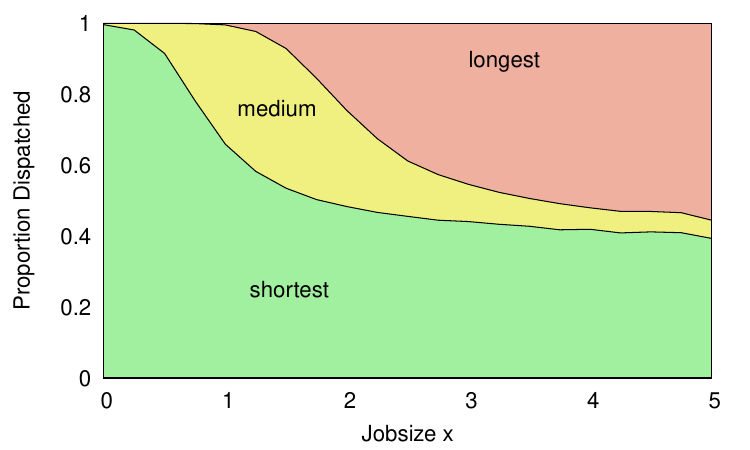}
  \end{tabular}
  \caption{Rank of the queue as a function of the job size $x$ with
    optimal policy for three servers and $\text{Exp}(1)$-distributed jobs.
    Short jobs are favored at the expense of longer jobs, i.e.\
    jobs are treated ``unfairly''.
  }
  \label{fig:qrank}
\end{figure}

Then let us consider which fraction of jobs of size $x$
are assigned to the shortest queue, to the queue with a medium length,
and to the longest queue.
This quantity can be seen as quality of service (QoS) metric.
Suppose we have $\sn=3$ servers, Poisson arrival process,
and Exp$(1)$-distributed jobs.
Figure~\ref{fig:qrank} depicts the resulting QoS metric
for $\rho=0.6$ and $\rho=0.9$.
When the load is low, almost all jobs, even long jobs,
are routed to the shortest queue.
(With $\rho=0.4$ this is practically the case and LWL is near optimal.)
However, as the load is high, $\rho=0.9$,
the situation has changed dramatically
and only the short jobs can enjoy the shortest queue.
For example, fewer than $50\%$ of jobs with size $x \ge 2$
are routed to the shortest queue.

\section{Conclusions}

This report provides technical improvements to the numerical solution methodology
developed for
the dispatching
problem in \cite{hyytia-peva-2024}, as well as
some new insights regarding the optimal dispatching rules itself.
Even though the optimality equations can be expressed in compact forms,
the resulting numerical problem remains computationally challenging.
In our examples, we had to limit ourselves to consider systems with a small number of servers,
say $\sn \le 6$.
The good news is that large systems behave fundamentally differently,
and near-optimal decisions can be based on finding an idle server, or a server that has a small backlog.
This is exactly what, e.g., the join-the-idle-queue (JIQ) does \cite{lu-peva-2011}.

\clearpage
\bibliographystyle{IEEEtran}
\bibliography{esa}

\end{document}